\documentclass{amsart}
\usepackage{amssymb,amsthm, amsmath}
\usepackage{bm}
\title{Some automorphism invariance properties for multicontractions}
\author{Chafiq Benhida} 
\address{UFR de Math\'ematiques, Universit\'e des Sciences et Technologies
de Lille, F-59655 Villeneuve D'Ascq Cedex, France}
\email{Chafiq.Benhida@math.univ-lille1.fr}

\author{Dan Timotin}
\address{Institute of Mathematics of the Romanian Academy, P.O. Box 1-764, Bucharest 
014700, Romania}
\email{Dan.Timotin@imar.ro}

\date{\today}

\newtheorem{theorem}{Theorem}[section]
\newtheorem{proposition}[theorem]{Proposition}
\newtheorem{lemma}[theorem]{Lemma}

\theoremstyle{definition}
\newtheorem{definition}[theorem]{Definition}

\newenvironment{tpmatrix}{\left(\begin{smallmatrix}}{\end{smallmatrix}\right)}

\let\mathbfit\bm

\def\J{{\mathbfit J}}

\def\L{{\mathbfit L}}
\def\M{{\mathbfit M}}

\def\DD{{\mathcal D}}
\def\EE{{\mathcal E}}

\def\GG{{\mathcal G}}
\def\HH{{\mathcal H}}

\def\KK{{\mathcal K}}
\def\LL{{\mathcal L}}
\def\MM{{\mathcal M}}
\def\NN{{\mathcal N}}

\def\UU{{\mathcal U}}

\def\XX{{\mathcal X}}
\def\YY{{\mathcal Y}}
\def\ZZ{{\mathcal Z}}

\def\Cfr{{\mathfrak C}}

\def\Ffr{{\mathfrak F}}

\def\Jfr{{\mathfrak J}}

\def\Lfr{{\mathfrak L}}

\def\Rfr{{\mathfrak R}}

\newcommand{\lfr}{\mathfrak l}
\newcommand{\rfr}{\mathfrak r}
\newcommand{\Flip}{F}

\def\BBB{{\mathbb B}}
\def\CCC{{\mathbb C}}

\def\FFF{{\mathbb F}}

\def\TTT{{\mathbb T}}

\def\0{{\mathbf 0}}
\def\1{{\mathbf 1}}

\def\<{\langle}
\def\>{\rangle}

\def\Om{Z}

\newcommand\Hilb{{\mathbf h}}
\newcommand\Fock{{\Ffr_n}}
\newcommand\Focks{{\Ffr_n^s}}

\newcommand\SOT{{\rm SOT}}
\newcommand\tFlip{{\tilde\Flip}}
\newcommand\Rr{r}

\begin{document}

\begin{abstract}
In the theory of row contractions on a Hilbert space, as initiated by Popescu, two important objects  are the Poisson kernel and the characteristic function. We determine their behaviour with respect to the action of the group of unitarily implemented automorphisms of the algebra generated by creation operators on the Fock space. The case of noncommutative varieties, introduced recently by Popescu, is also discussed.
\end{abstract}

\maketitle

\numberwithin{equation}{section}

\section{Introduction}\label{se:intro}

Among the attempts to
extend the dilation theory of contractions on a Hilbert space, as developed
in~\cite{SNF}, to multivariable operator theory, a most notable achievement is the theory of row contractions (which we will call below \emph{multicontractions}),
initiated by the 
work of Gelu Popescu~\cite{GP1, GP, GP1.5}. It  has been pursued
in the last two decades by Popescu and others (see for instance,~\cite{AP, AP2, DP, DP2, GP2, GP2.5}).
Popescu's theory is essentially noncommutative; later, starting with~\cite{Ar}, interest
has developed around the case of commuting multioperators. This presents some specific features; on the other hand, many properties of the commuting case can be obtained
from the noncommuting situation.
The two recent papers of Popescu~\cite{GP3, GP4} have pursued systematically
the development of the commutative situation from the noncommutative one, 
putting it into the more general framework of \emph{constrained multioperators}
(see Section~\ref{se:constrained}). 

Two objects related to a multicontraction play a significant role in Popescu's 
theory: the \emph{Poisson kernel} and the \emph{characteristic function}. The Poisson kernel
is an important tool used by Popescu in order to prove the von Neumann inequality
for row contractions; the characteristic function is essential
in the model theory of completely noncoisometric multicontractions, while
its commutative counterpart is related to 
Arveson's curvature~\cite{Ar, BES, GP3}. 

On the other hand, as the Sz.-Nagy--Foias theory of single contractions is related
to classical function theory in the unit circle, an analoguous role
in Popescu's theory is played by algebras generated by creation operators
on the Fock space. There is a distinguished group of automorphisms of these
algebras, that have been introduced by Voiculescu in~\cite{V} (and discussed
recently in~\cite{DP}); they
are the noncommutative analogues of the analytic automorphisms of the
unit ball. These automorphisms act on multicontractions, and it is interesting
to see what is the effect of this transformations on the Poisson kernel
and on the characteristic function. The purpose of this paper is
to show, firstly, that these objects obey natural rules of transformation,
and secondly, that the rules of transformation also extend to the case
of constrained multicontractions. In particular, the relation between
the transformation rules for commuting and for noncommuting
multicontractions is clarified.

The first three sections following the 
introduction contain mostly preliminary material. The main results,
Theorems~\ref{th:theta} and~\ref{th:theta2}, are proved in Sections~\ref{se:ma}
and~\ref{se:constrained}. In connection to constrained objects, the last section discusses invariant ideals of the noncommutative Toeplitz algebra.

The method of proof uses, in order to avoid complicated computations,
the machinery of Redheffer products. This machinery may seem
unfamiliar, but it provides a simple and short way to reach
the main results. Up to a certain point, using Redheffer
products is equivalent to composing $J$-unitary operators,
but there is a slightly larger level of generality that 
happens to be important in our context. For an illuminating discussion
of these facts, see the first two sections of~\cite{Y}.

\section{Preliminaries}\label{se:char}

\subsection{The Fock space}
A main object of study is formed by the Fock space and the non-commutative Toeplitz
algebras that act on it.
We will follow mainly the work of Popescu~\cite{GP1, GP, GP2}, as well as~\cite{DP}.

In the whole paper we will fix a positive integer $n$.
We denote by $\FFF_n^+$  the free semigroup with the $n$ generators $1, \dots, n$ and unit $\emptyset$.
An element $w=i_1\cdots i_k\in \FFF_n^+$  is  called a \emph{word} in the letters $1, \dots, n$,
and its \emph{length} is $|w|=k$. 

If $A=(A_1, \dots, A_n)\in\LL(\HH)^n$ is a (not necessarily commuting) multioperator, 
we denote $A_\emptyset=I_\HH$ and, if,
$w=i_1\cdots i_k\in \FFF_n^+$, then 
$A_w=A_{i_1}\cdots A_{i_k}\in\LL(\HH)$.

Consider
an $n$-dimensional complex Hilbert space $\Hilb_n$, with basis vectors
$e_1,\dots, e_n$. The \emph{full Fock space} is then
\[
\Fock=\bigoplus_{k\ge 0} \Hilb_n^{\otimes k}
\]
where $\Hilb_n^{\otimes 0}=\CCC\1$ and $\Hilb_n^{\otimes k}$ is the tensor product
of $k$ copies of $\Hilb_n$. 
An orthonormal basis of $\Fock$ is given by $(e_w)_{w\in\FFF_n^+}$, where
$e_\emptyset=\1\in \Hilb_n^{\otimes 0}$, while, if $w=i_1\cdots i_k\in \FFF_n^+$,
then $e_w= e_{i_1}\otimes\cdots\otimes e_{i_k}\in \Hilb_n^{\otimes k}$.

The left creation operators $L_i\in\LL(\Fock)$,
$i=1,\dots n$, are defined by 
\[
L_i\xi= e_i\otimes \xi, \quad \xi\in \Fock.
\]
The norm closed algebra generated by $L_1, \dots, L_n$ is denoted by $\lfr_n$,
and the weakly closed algebra by $\Lfr_n$.

Similarly, we have right creation operators $R_i$ given by 
\[
R_i\xi=\xi\otimes e_i,
\]
while the norm closed and weakly closed algebras they generate are denoted by
$\rfr_n$ and $\Rfr_n$ respectively. Each of the algebras $\Lfr_n$ and $\Rfr_n$ is
the commutant of the other.

We
can write any $f\in\Lfr_n$ as a formal series
$f=\sum_w \hat f_w L_w$. For each
$r<1$ the series 
\begin{equation}\label{eq:f_r}
f_r:=\sum_w \hat f_w r^{|w|}L_w
\end{equation}
converges uniformly, and thus $f_r\in\lfr_n$; then $f=\SOT-\lim_{r\to1}f_r$.
A similar statement is valid for $\Rfr_n$.

The \emph{flip} operator is
the involutive unitary $\Flip\in\LL(\Fock)$ which acts on simple tensors by reversing
the order of the components:
\[
\Flip(e_{i_1}\otimes e_{i_2}\otimes \cdots \otimes e_{i_n})
= e_{i_n}\otimes \cdots \otimes e_{i_2}\otimes e_{i_1}.
\]
We have then $R_i=\Flip L_i\Flip$.

If
$\EE, \EE_*$ are Hilbert spaces, then a linear operator $M:\Fock\otimes\EE\to\Fock\otimes\EE_*$
is called \emph{multianalytic} if
\[
M(L_i\otimes I_\EE)=(L_i\otimes I_{\EE_*})M \quad \forall i=1,\dots,n.
\]
$M$ is then uniquely determined by the ``coefficients" $m_w\in\LL(\EE, \EE_*)$,
defined by
\begin{equation*}
\<m_{\tilde w} k,k'\>= \<M(\1\otimes k), e_w\otimes k_*\>, \quad
k\in\EE, k_*\in\EE_*, w\in \FFF_n^+,
\end{equation*}
where $\tilde w$ is the reverse of $w$, i.e., $\tilde w =i_k\cdots i_1$ if $ w=i_1\cdots i_k$; we can then
associate with $M$ the formal Fourier expansion
\begin{equation}\label{eq:formal}
\hat M(R_1, \dots , R_n)=\sum_{w\in\FFF_n^+} R_w\otimes m_w.
\end{equation}
%


%
%
%
%

\subsection{Redheffer products}\label{se:Redheffer}

Several computations that appear in the sequel can be gathered
in a simple uniform framework if we use the formalism of Redheffer products.
The basic reference
is~\cite{Re}; we will follow the exposition in~\cite{T}.

Suppose 
\[
\L=\begin{pmatrix} A & B\\ C & D
\end{pmatrix}, \quad 
\L_1=\begin{pmatrix} A_1 & B_1\\ C_1 & D_1
\end{pmatrix}
\]
are bounded operators mapping $\XX\oplus \UU$ (respectively 
$\XX_1\oplus \UU_1$) into $\YY\oplus \ZZ$ (respectively 
$\YY_1\oplus \ZZ_1$); also $\UU_1=\ZZ$ and $\XX=\YY_1$. 
Under the  assumption 
\[
I-B_1C \text{ is invertible}\tag{*}
\]
it follows that $I-CB_1$ is also invertible, and we define the \emph{Redheffer product} 
by
\begin{equation}\label{eq:rprod}
\M=\L\circ \L_1=
\begin{pmatrix}
A(I-B_1C)^{-1}A_1 & B+A(I-B_1C)^{-1}B_1D\\
C_1+D_1C(I-B_1C)^{-1}A_1 & D_1(I-CB_1)^{-1}D
\end{pmatrix}.
\end{equation}
$\M$ is an operator from $\XX_1\oplus \UU$ to $\YY\oplus \ZZ_1$. It is useful
tu visualize the interlacing of spaces by input-output boxes, in a manner
suggested by system theory:

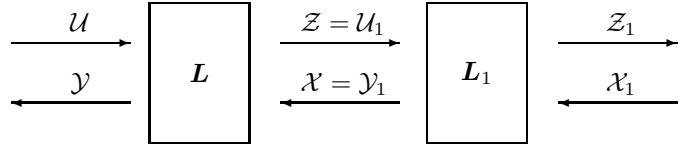
\setlength{\unitlength}{1.5pt}
\newcommand{\Mybox}{\framebox(25,35)}
\begin{center}
\begin{picture}(100,50)(50,10)

\put(10,40){\vector(1,0){30}}
\put(40,25){\vector(-1,0){30}}
\put(45,15){\Mybox{$\L$}}

\put(25,43){$\UU$}
\put(25,28){$\YY$}

\put(78,40){\vector(1,0){30}}
\put(108,25){\vector(-1,0){30}}
\put(115,15){\Mybox{$\L_1$}}

\put(83,43){$\ZZ=\UU_1$}
\put(83,28){$\XX=\YY_1$}

\put(148,40){\vector(1,0){30}}
\put(178,25){\vector(-1,0){30}}

\put(160,43){$\ZZ_1$}
\put(160,28){$\XX_1$}

\end{picture}
\end{center}

We will write also $\beta_\L(A_1, B_1)$ and $\alpha_\L(B_1)$ for the entries in
the first row of $\L\circ \L_1$ (as given by~\eqref{eq:rprod}).

The basic properties of the Redheffer product are gathered
in the following proposition. In its statement it is tacitly assumed
that condition~(*) is satisfied, when necessary.

\begin{proposition}\label{pr:red}
(i) The identities matrices (on the corresponding spaces) act as unit 
elements also for the Redheffer products.

(ii) If $\L$ is invertible, $\L_1=\L^{-1}$, and one can form $\L\circ \L_1$, then
$\L_1$ is also the inverse of $\L$ with respect to the Redheffer product.

(iii) The Redheffer product is associative: if $\L,\L_1, \L_2$ are given, and all Redheffer
products in~\eqref{eq:assoc} can be formed, then
\begin{equation}\label{eq:assoc}
\L\circ (\L_1\circ \L_2)= (\L\circ \L_1)\circ \L_2.
\end{equation}

(iv) $\L, \L_1$ contractions (isometries, coisometries, unitaries) imply
$\L\circ \L_1$ contraction (isometry, coisometry, unitary respectively). In particular,
if $\L$ and $B_1$ are contractions, then $\alpha_\L(B_1)$ is also a contraction.
\end{proposition}

A particular case that will be useful is $\ZZ_1=\{0\}$ (and thus $C_1=D_1=0$).

In connection with Redheffer products, we need also a lemma
concerning the structure of unitary $2\times 2$ matrices. To state it, remember that if
$\EE_1, \EE_2$ are two Hilbert spaces, and $C:\EE_1\to \EE_2$ is a contraction,
one defines the \emph{defect operator} $D_C=(\1_{\EE_1}-C^*C)^{1/2} \in\LL(\EE_1)$ and the
\emph{defect space} $\DD_C=\overline{D_C \EE_1}\subset \EE_1$.

\begin{lemma}\label{le:structure}
A $2\times 2$ operator matrix from $\EE_1\oplus\EE_2$
to $\EE'_1\oplus\EE'_2$    that has $A^*$ as its $(2,1)$ entry, while the $(1,1)$ 
entry has dense range, has the form
\begin{equation}\label{eq:structure}
\J=\begin{pmatrix}
Z_*D_{A^*} & -Z_*AZ^*\\
A^* & D_A Z^*
\end{pmatrix},
\end{equation}
$Z_*:\DD_{A^*}\to \EE'_1$ and $Z:\DD_A\to \EE_2$ being unitary operators.
\end{lemma}

\section{Automorphisms}\label{se:autoball}

The analytic automorphisms of the unit ball $\BBB^n$ act by composition on any 
Hilbert space of functions on $\BBB^n$. There exist corresponding
unitarily implemented automorphisms on the non-commutative 
Toeplitz algebras on the Fock space. 

\subsection{The commutative case: automorphisms of the unit ball}\label{sse:autoball}

There are two different descriptions of the automorphisms of the 
unit ball $\BBB^n$. In view of
further extensions, we  identify elements 
in  $\BBB^n$ with row $1\times n$ contractive matrices.
Naturally, the action of an $n\times n$ matrix on such
an element will be done by multiplication on the right.

\textit{First form}. Start with the group 
$U(1,n)$ of $(n+1)\times (n+1)$ matrices $X$ that are
$J$-unitary,
where $J=\begin{pmatrix} -1 & 0\\0 & I_n\end{pmatrix}$;
that is, $X^*JX=J$. 
According to the decomposition 
$\CCC^{n+1}=\CCC\oplus \CCC^n$, one writes $X=\begin{pmatrix} x & {y}\\ {z}^t& X'\end{pmatrix}$; note that with these conventions $x$ is a scalar, while ${y}$  and ${z}$ are 
row matrices.
Accordingly, there is a corresponding map $\phi_X:\BBB^n\to \BBB^n$,
defined by 
\begin{equation}\label{eq:phiX}
\phi_X(\lambda)=(x-\lambda{z}^t)^{-1}(\lambda X'-{y}).
\end{equation}

Then the map $X\mapsto \phi_X$ is a group antihomomorphism from
$U(1,n)$ to the group af automorphisms of $\BBB^n$ (the ``anti"
being due to our decision to see elements of $\BBB^n$ as row matrices and
write the action of the group on the right);
this antihomorphism is onto, and its kernel is formed by scalar
unitaries.

\textit{Second form}. A variant of~\eqref{eq:phiX} which uses
a unitary instead of a $J$-unitary  matrix is more natural
in the context of Redheffer products. Namely,
if $Y=\left(\begin{smallmatrix}
a & b\\ c& d
\end{smallmatrix}\right)$ is a unitary $(n+1)\times(n+1)$ matrix,
then we can consider the map (see~\eqref{eq:rprod})
\begin{equation}\label{eq:phiY}
\alpha_Y(\lambda)=b+a\lambda (I-c\lambda)^{-1}d.
\end{equation}
The corresponding diagram is

\setlength{\unitlength}{1.5pt}
\begin{center}
\begin{picture}(100,50)(50,10)

\put(10,40){\vector(1,0){30}}
\put(40,25){\vector(-1,0){30}}
\put(45,15){\Mybox{$\L_Y$}}

\put(25,43){$\CCC^n$}
\put(25,28){$\CCC$}

\put(78,40){\vector(1,0){30}}
\put(108,25){\vector(-1,0){30}}
\put(115,15){\Mybox{$\L_1$}}

\put(93,43){$\CCC^n$}
\put(93,28){$\CCC$}

\put(148,40){\vector(1,0){30}}
\put(178,25){\vector(-1,0){30}}

\put(160,43){$\{0\}$}
\put(160,28){$\CCC$}

\end{picture}
\end{center}
where $\L_1=\left(\begin{smallmatrix} 1 & \lambda\\ 0&0
\end{smallmatrix}\right):\CCC\oplus \CCC^n\to\CCC\oplus \{0\}$.

By Proposition~\ref{pr:red} (iv), for $\lambda$ contractive, $\alpha_Y(\lambda)$
is also contractive. Thus $\alpha_Y$ 
is an analytic map from $\BBB^n$ to $\BBB^n$; it is even an automorphism,
since again Proposition~\ref{pr:red} (ii) implies
$\alpha_Y^{-1}=\alpha_{Y^*}$.

The passage from~\eqref{eq:phiX} to~\eqref{eq:phiY} is done by the
formulas
\begin{equation}\label{eq:relXY}
a=x^{-1}, \qquad b=-x^{-1}{y},
\qquad c=x^{-1}{z}^t,\qquad d= X'- x^{-1}{z}^t {y}.
\end{equation}
These formulas can be inverted, provided $a\not=0$.

Working in the context of Redheffer products,~\eqref{eq:phiY} is 
more convenient; however,~\eqref{eq:phiX} is related to the automorphisms
in Subsection~\ref{sse:implemented}. Also, while in~\eqref{eq:phiX} any $J$-unitary
produces an automorphism, in~\eqref{eq:phiY} we must require $a\not=0$.

\subsection{The noncommutative case: the Fock space}\label{sse:implemented}

We shall introduce some facts and notations from~\cite{DP}; in
Section~4 therein
the automorphisms of the algebra $\Lfr_n$ are investigated.
It is shown that all contractive automorphisms of $\Lfr_n$ are
actually unitarily implemented, and they are also automorphisms
of the $C^*$-algebra~$\lfr_n$. 

A detailed
description of these automorphisms can be obtained following~\cite{V}. 
As in Section~\ref{se:autoball}, take $X\in U(1,n)$, 
$X=\begin{tpmatrix} x & {y}\\ {z}^t & X'\end{tpmatrix}$.
Write also $L[\zeta]=\sum_{i=1}^n \zeta_iL_i$ for $\zeta\in\CCC^n$. Then
 there is an automorphism $\Phi_X$ of $\Lfr_n$
such that the restriction to the generators is given by
\begin{equation}\label{eq:voic1}
\Phi_X (L[\zeta])=(xI-L[{z}])^{-1}(L[X'\zeta]-(\zeta\cdot {y}^t) I).
\end{equation}
This automorphism is implemented by a unitary $U_X\in\LL(\Fock)$, which
satisfies
\begin{equation}\label{eq:voic2}
U_X(A{e_\emptyset})=\Phi_X(A) (xI-L[z])^{-1}{e_\emptyset}
\end{equation}
for all $A\in\Lfr_n$; this means that $\Phi_X(A)=U_XAU_X^*$ for all $A\in\Lfr_n$.
The map $X\mapsto \Phi_X$ from $U(1,n)$ to the automorphisms of $\Lfr_n$ 
has as image all unitarily implemented automorphisms (which actually
coincide with all
contractive automorphisms), and its kernel consists of the scalar
matrices $xI_{n+1}$, with $x\in\TTT$.

To make the connection with~\ref{sse:autoball}, apply~\eqref{eq:voic1} 
for $\zeta$ a basis vector; one obtains
\[
\Phi_X (L_i)=(x I-L[{z}])^{-1}(\sum_{j=1}^n X'_{ji} L_j-{y}_{i} I),
\]
while writing~\eqref{eq:phiX} on coordinates yields
\[
(\phi_X(\lambda))_i=(x-\sum_{j=1}^n {z}_j\lambda_j)^{-1}
(\sum_{j=1}^n \lambda_j X'_{ji}-{y}_i).
\]
Consequently,~\eqref{eq:voic1} can be obtained by formally replacing
$\lambda_i$ in~\eqref{eq:phiX} with $L_i$.

One can interpret also these automorphisms in terms of Redheffer products.
Suppose that, as in~\ref{sse:autoball}, on defines the unitary matrix
$Y=\begin{tpmatrix}
a & b\\ c& d
\end{tpmatrix}\in\LL(\CCC\oplus \CCC^n)$ by~\eqref{eq:relXY};
denote by $\iota:\CCC\to\Fock$ the inclusion map that sends
$1$ to ${e_\emptyset}$,
and
\[
\begin{split}
&\L_Y=
\begin{pmatrix}
a I_{\Fock }  & b\otimes I_{\Fock }\\
c\otimes I_{\Fock}& d\otimes I_{\Fock}
\end{pmatrix}: 
\Fock\oplus (\Fock\otimes \CCC^n)
\to
\Fock \oplus (\Fock\otimes \CCC^n)
\\
&\L_1=
\begin{pmatrix}
\iota  & L \\
0 & 0
\end{pmatrix}: \CCC\oplus (\Fock\otimes \CCC^n)
\to 
(\Fock)\oplus \{0\} 
\end{split}
\]
(we have implicitely used the fact that we can identify
$\Fock\otimes \CCC^n$ with $\Ffr_n^n$).

\setlength{\unitlength}{1.5pt}
\begin{center}
\begin{picture}(100,50)(50,10)

\put(10,40){\vector(1,0){30}}
\put(40,25){\vector(-1,0){30}}
\put(45,15){\Mybox{$\L_Y$}}

\put(15,43){$\Fock\otimes\CCC^n$}
\put(22,28){$\Fock$}

\put(78,40){\vector(1,0){30}}
\put(108,25){\vector(-1,0){30}}
\put(115,15){\Mybox{$\L_1$}}

\put(83,43){$\Fock\otimes \CCC^n$}
\put(90,28){$\Fock$}

\put(148,40){\vector(1,0){30}}
\put(178,25){\vector(-1,0){30}}

\put(160,43){$\{0\}$}
\put(160,28){$\CCC$}

\end{picture}
\end{center}

Then it follows immediately from the discussion above that:
\begin{equation}\label{eq:alphi}
\alpha_{\L_Y}(L)= \Phi_X(L)= U_XLU_X^*.
\end{equation}
Moreover, from~\eqref{eq:voic2} and~\eqref{eq:relXY} we have
\[
U_X({e_\emptyset})= x^{-1}(I-x^{-1}\sum_{j=1}^n {z}_jL_j)^{-1}{e_\emptyset}
=a(I-\sum_{j=1}^n c_j L_j)^{-1}{e_\emptyset},
\]
whence
\[
\beta_{\L_Y}(\iota, L)= U_X\iota.
\]
Let us also note that, if $X$ and $Y$ are related by~\eqref{eq:relXY}, then 
\begin{equation}
	\Phi_X^{-1}(L)=\Phi_{X^{-1}}(L)=\alpha_{\L_{Y^*}}(L).
	\label{eq:inverse}
\end{equation}


\smallskip
Suppose now that we want to obtain similar relations with $R$ instead of $L$. We
may immediately note that $R_i=\Flip L_i\Flip$, which  leads
to
\[
\alpha_{\L_Y}(R)=\Flip U_X \Flip R \Flip U_X^* \Flip.
\]
But we can actually say more. 
Since $\Rfr_n$ is the commutant of $\Lfr_n$, it follows
that $U_XBU_X^*\in\Rfr_n$ for all $B\in\Rfr_n$; this can
be made precise using the following lemma.

\begin{lemma}\label{le:uflip}
For all $X\in U(1,n)$, we have $U_X\Flip=\Flip U_X$.
\end{lemma}

\begin{proof}
We have to check the relation on simple tensors $e_w$, where $w=i_1\cdots i_k$.
Denote also $\tilde w=i_k\cdots i_1$.
According to~\eqref{eq:voic2}, we have
\[
\begin{split}
U_X\Flip(e_w)&=U_X(e_{\tilde w})=U_X(L_{i_k}\cdots L_{i_1}{e_\emptyset})\\
&=\Phi_X(L_{i_k}\cdots L_{i_1}) (xI-L[{z}])^{-1}{e_\emptyset}\\
&=\Phi_X(L_{i_k})\cdots\Phi_X( L_{i_1}) (xI-L[{z}])^{-1}{e_\emptyset}.
\end{split}
\]
If we write $\tfrac{1}{x} {z}=(a_1,\dots, a_n)$, we have
\[
(xI-L_{{z}})^{-1}=x\sum_{\text{all words }j=j_1\dots j_s}a_{j_1}L_{j_1}\dots a_{j_s} L_{j_s},
\]
where the sum is norm convergent.

To simplify notations, define the map $\tFlip:\Lfr_n\to\Lfr_n$ 
by the formula $\tFlip (L_v)=L_{\tilde v}$. Then:
\begin{enumerate}
\item $\tFlip(AB)=\tFlip(B)\tFlip(A)$;
\item $\tFlip(L_v)=L_v$ if $v$ has length 0 or 1;
\item
$\tFlip ((xI-L[{z}])^{-1})=(xI-L[{z}])^{-1}$;
\item $\Flip L_v {e_\emptyset}=\tilde\Flip(L_v){e_\emptyset}$.
\end{enumerate}

Therefore, applying~\eqref{eq:voic1}, we have
\[
\begin{split}
&\Phi_X(L_{i_k})\cdots\Phi_X( L_{i_1}) (xI-L[{z}])^{-1}{e_\emptyset}\\
&\qquad= (xI-L[{z}])^{-1} (L[X'\zeta_{i_k}]-\< \zeta_{i_k}, {y}\> I)
(xI-L[{z}])^{-1}\dots \\
&\qquad\qquad
(xI-L[{z}])^{-1}
(L[X'\zeta_{i_1}]-\< \zeta_{i_1}, {y}\> I)
(xI-L[{z}])^{-1}{e_\emptyset}\\
&\qquad= \tilde\Flip\big[(xI-L[{z}])^{-1} (L[X'\zeta_{i_1}]-\< \zeta_{i_1}, {y}\> I)
(xI-L[{z}])^{-1}\dots \\
&\qquad\qquad
(xI-L[{z}])^{-1}
(L[X'\zeta_{i_k}]-\< \zeta_{i_k}, {y}\> I)
(xI-L[{z}])^{-1}\big]{e_\emptyset}\\
&\qquad= \Flip(xI-L[{z}])^{-1} (L[X'\zeta_{i_1}]-\< \zeta_{i_1}, {y}\> I)
(xI-L[{z}])^{-1}\dots \\
&\qquad\qquad
(xI-L[{z}])^{-1}
(L[X'\zeta_{i_k}]-\< \zeta_{i_k}, {y}\> I)
(xI-L[{z}])^{-1}{e_\emptyset}\\
&\qquad =\Flip \Phi_X(L_{i_1})\cdots\Phi_X( L_{i_k}) (xI-L[{z}])^{-1}{e_\emptyset}\\
&\qquad=\Flip U_X(e_w).
\end{split}
\]
The lemma is proved.
\end{proof}

As a consequence, $\Flip U_X\Flip= U_X$, and we have 
\[
\alpha_{\L_Y}(R)=U_X R U_X^*.
\]


\section{Multicontractions}\label{se:multi}

Suppose $T=(T_1,\dots, T_n)\in\LL(\HH)^n$ is a  multicontraction;
that is, 
\[
\sum_{i=1}^n T_iT_i^*\le 1_\HH. 
\]
This is the same as requiring 
the row operator $T=(T_1\ \cdots\ T_n):\HH^n\to\HH$ to be a contraction. (We will
currently denote with the same letter $T$ the multioperator and the
associated row contraction.) 
Accordingly, we have the operators
$D_T=(\1_{\HH^n}-T^*T)^{1/2}$ and $D_{T^*}=(\1_\HH-TT^*)^{1/2}$, 
and the spaces $\DD_T=\overline{D_T\HH^n}\subset \HH^n$, 
$\DD_{T^*}=\overline{D_{T^*}\HH}\subset\HH$. 
If the row operator $T$ is a strict contraction, we will say that $T$ is a \emph{strict
multicontraction}.

There exists an $\lfr_n$-functional calculus for a multicontraction $T$;
it is the unique completely contractive homomorphism $\rho:\lfr_n\to \LL(\HH)$,
such that $\rho(L_i)=T_i$. This homomorphism can be extended to $\Lfr_n$ in 
an important particular case. Namely, $T$ is called \emph{completely noncoisometric (c.n.c.)}
if there is no $h\in\HH$, $h\not=0$, such that
\[
\sum_{|w|=k}\|T^*_w h\|^2=\|h\|^2\quad\mbox{for all $k\ge0$}. 
\]
If $T$ is c.n.c., then $\rho$ can be extended to a completely contractive homomorphism defined 
on $\Lfr_n$, that
we will denote with the same letter, $\rho:\Lfr_n\to \LL(\HH)$~\cite{GP2}. 
If $f\in\Lfr_n$, then $f_r\in\lfr_n$, and we may apply $\rho$ to obtain 
\[
\rho(f_r)=\sum_w \hat f_w r^{|w|}T_w
\]
with the sum on the right converging absolutely. If $T$ is c.n.c.,
then we have also
\[
\rho(f)=\SOT-\lim_{r\to1} \rho(f_r).
\]

Similar results are valid for $\rfr_n$ and $\Rfr_n$, the corresponding
functional calculus being denoted by $\rho'$.

The next definition introduces two basic objects that appear in Popescu's theory of multicontractions (see \cite{GP2.5, GP3}).

\begin{definition}\label{de:poissontheta}
Suppose $T$ is a multicontraction. Then:

(a) The \emph{Poisson kernel} $K_T$ is the operator
$K_T:\HH\to\Fock\otimes \DD_{T^*}$
defined by
\[
K_T h =\sum_w e_w \otimes D_{T^*} T_w^*h.
\]

(b) The \emph{characteristic function} 
 $\Theta_T$ is the multianalytic operator
\[
\Theta_T: \Fock\otimes \DD_T\to\Fock\otimes \DD_{T^*}
\]
having the formal Fourier representation
\begin{equation}\label{eq:theta1}
\begin{split}
\hat\Theta_T(R_1,\cdots,R_n)&=-I_{\Fock}\otimes T +(I_{\Fock}\otimes D_{T^*})
\left( I_{\Fock\otimes\HH} -\sum_{i=1}^n R_i\otimes T_i^*\right)^{-1}\\
 &\qquad\qquad[R_1\otimes I_\HH, \dots , R_n\otimes I_\HH] (I_{\Fock}\otimes D_T)
\big|\Fock\otimes\DD_T.
\end{split}
\end{equation}
\end{definition}

The following proposition gathers several results from~\cite{GP3}.

\begin{proposition}\label{pr:gp3} 
{\rm(i)}  The Poisson kernel and the characteristic function are contractions, and
$K_T K_T^*+\Theta_T \Theta_T^*=I_{\Fock\otimes \DD_{T^*}}$. 

{\rm (ii)} If we define, for $0<r\le1$,
\[
K_{T,r} h =\sum_w r^{|w|} e_w \otimes D_{T^*} T_w^*h
=(I_{\Fock}\otimes D_{T^*})
\left( I_{\Fock\otimes\HH} -\sum_{i=1}^n r R_i\otimes T_i^*\right)^{-1}(e_\emptyset\otimes h),
\]
then 
\begin{equation}\label{eq:Kr}
K_T=\SOT-\lim_{r\to1} K_{T,r}.
\end{equation}

{\rm (iii)}
If we replace $R_i$ with $rR_i$, $0<r<1$, in~\eqref{eq:theta1}, then the inverse in the
right hand side exists,  the equation can be used to define 
$\Theta_T(rR)$, and 
\begin{equation}\label{eq:Thr}
\Theta_T=\SOT-\lim_{r\to1} \Theta_T(rR).
\end{equation}
\end{proposition}

We can interpret the Poisson kernel and the characteristic function
by means of the Redheffer product (see Section~\ref{se:Redheffer}). 
Remember that $\iota:\CCC\to\Fock$ is the embedding $z\mapsto z e_\emptyset$.
Take $r<1$, and define then
\begin{align}\label{eq:defLT}
\L_T=
\begin{pmatrix}
I_\Fock \otimes D_{T^*}  & -I_\Fock \otimes T\\
I_\Fock\otimes T^* & I_\Fock\otimes D_T
\end{pmatrix}: 
(\Fock\otimes \HH)\oplus (\Fock\otimes \DD_T)
&\to
(\Fock\otimes \DD_{T^*}) \oplus (\Fock\otimes \HH^n)
\\
\L_\Rr=
\begin{pmatrix}
\iota\otimes I_\HH  & rR \otimes I_\HH\\
0 & 0
\end{pmatrix}: \HH\oplus (\Fock\otimes \HH^n)
&\to 
(\Fock\otimes \HH)\oplus \{0\} \label{eq:defLr}
\end{align}

\setlength{\unitlength}{1.5pt}
\begin{center}
\begin{picture}(100,50)(50,10)

\put(10,40){\vector(1,0){30}}
\put(40,25){\vector(-1,0){30}}
\put(45,15){\Mybox{$\L_T$}}

\put(15,43){$\Fock\otimes\DD_T$}
\put(15,28){$\Fock\otimes\DD_{T^*}$}

\put(78,40){\vector(1,0){30}}
\put(108,25){\vector(-1,0){30}}
\put(115,15){\Mybox{$\L_\Rr$}}

\put(83,43){$\Fock\otimes \HH^n$}
\put(83,28){$\Fock\otimes\HH$}

\put(148,40){\vector(1,0){30}}
\put(178,25){\vector(-1,0){30}}

\put(160,43){$\{0\}$}
\put(160,28){$\HH$}

\end{picture}
\end{center}

Then from~\eqref{eq:rprod} it follows that
\[
\L_T\circ \L_\Rr=
\begin{pmatrix}
K_{T,r} & \Theta_T(rR) \\ 0 & 0 
\end{pmatrix}
:\HH\oplus (\Fock\otimes \DD_T)\to (\Fock\otimes \DD_{T^*})\oplus \{0\}.
\]
Otherwise stated,
\begin{equation}\label{eq:KTh}
\Theta_T(rR)=\alpha_{\L_T}(rR \otimes I_\HH), 
\qquad K_{T,r}=\beta_{\L_T} (\iota\otimes I_\HH, rR \otimes I_\HH).
\end{equation}

%
%
%

%

%

\section{Multicontractions and automorphisms}\label{se:ma}

By using the functional calculus $\rho$ (see Section~\ref{se:multi}), we can extend the action of automorphisms $\Phi_X$ to a multicontraction
$T$. This is done by defining
\[
\Phi_X(T)=\rho(\Phi_X(L)),
\]
and it follows from~\eqref{eq:alphi} that we have then also
\[
\Phi_X(T)=\alpha_{\L_Y}(T),
\]
where, as usually, $Y$ is connected to $X$ by formulas~\eqref{eq:relXY}.
Since the functional calculus $\rho$ is completely contractive, $\Phi_X(T)$
is also a multicontraction. According to~\eqref{eq:inverse}, we have also
\begin{equation}
	\Phi_X^{-1}(T)=\alpha_{\L_{Y^*}}(T).
	\label{eq:inverse2}
\end{equation}

The main result of this section is given by the next theorem.

\begin{theorem}\label{th:theta} For each $X$ there exist unitary operators 
$\Om:\DD_{\Phi_X^{-1}(T)}\to \DD_T$ and 
$\Om_*:\DD_{\Phi_X^{-1}(T)^*}\to \DD_{T^*}$, such 
that:

(i) $\Theta_{ \Phi_X^{-1}(T)}=( U_X \otimes \Om_*^* )\Theta_T ( U_X^*  \otimes \Om  )$.

(ii) $K_{\Phi_X^{-1} (T)} =( U_X \otimes \Om_*^* ) K_T$.
\end{theorem}

\begin{proof}
Let us define $\L_T$ and $\L_\Rr$ by formulas~\eqref{eq:defLT} and~\eqref{eq:defLr}
respectively, and $\L_X$ by
\[
\L_X=
\begin{pmatrix}
a \otimes I_{\Fock\otimes  \HH }  & b\otimes I_{\Fock \otimes \HH }\\
c\otimes I_{\Fock\otimes \HH}& d\otimes I_{\Fock\otimes \HH}
\end{pmatrix}: 
(\Fock\otimes \HH)\oplus (\Fock\otimes \HH^n)
\to
(\Fock\otimes \HH) \oplus (\Fock\otimes \HH^n)
\]
where $a,b,c,d$ are related to $X$ by formulas~\eqref{eq:relXY}. 

\setlength{\unitlength}{1.5pt}
\begin{center}
\begin{picture}(100,50)(80,10)

\put(10,40){\vector(1,0){30}}
\put(40,25){\vector(-1,0){30}}
\put(45,15){\Mybox{$\L_T$}}

\put(15,43){$\Fock\otimes\DD_T$}
\put(15,28){$\Fock\otimes\DD_{T^*}$}

\put(78,40){\vector(1,0){30}}
\put(108,25){\vector(-1,0){30}}
\put(115,15){\Mybox{$\L_X$}}

\put(83,43){$\Fock\otimes \HH^n$}
\put(83,28){$\Fock\otimes\HH$}

\put(148,40){\vector(1,0){30}}
\put(178,25){\vector(-1,0){30}}

\put(150,43){$\Fock\otimes \HH^n$}
\put(150,28){$\Fock\otimes\HH$}

\put(185,15){\Mybox{$\L_\Rr$}}

\put(215,40){\vector(1,0){30}}
\put(245,25){\vector(-1,0){30}}

\put(225,43){$\{0\}$}
\put(225,28){$\HH$}

\end{picture}
\end{center}

We want to apply the associativity of the Redheffer product, as stated in 
Proposition~\ref{pr:red} (iii):
\begin{equation}\label{eq:assoc2}
(\L_T\circ \L_X)\circ \L_\Rr= \L_T\circ (\L_X\circ \L_\Rr).
\end{equation}

First, we have
\[
\L_X\circ \L_\Rr =
\begin{pmatrix}
(U_X\iota) \otimes I_\HH & (U_XrRU_X^*)\otimes I_\HH\\
0 & 0
\end{pmatrix}.
\]
It follows that
\begin{equation*}
\begin{split}
\alpha_{\L_T}((U_XRU_X^*)\otimes I_\HH)&=-I_{\Fock}\otimes T +(I_{\Fock}\otimes D_{T^*})
\left( I_{\Fock\otimes\HH} -\sum_{i=1}^n (U_XrR_iU_X^*)\otimes T_i^*\right)^{-1}\\
 &\quad[U_X rR_1 U_X^*\otimes I_\HH, \dots ,U_X r R_n U_X^*\otimes I_\HH] (I_{\Fock}\otimes D_T)
\big|\Fock\otimes\DD_T\\
&= (U_X\otimes I_{\DD_{T^*}}) \Theta_T(rR) (U_X^*\otimes I_{\DD_T})
\end{split}
\end{equation*}
and
\begin{equation*}
\begin{split}
&\beta_{\L_T}\big((U_X\iota) \otimes I_\HH , (U_XrRU_X^*)\otimes I_\HH\big)\\
&\qquad=(I_{\Fock}\otimes D_{T^*})
\left( I_{\Fock\otimes\HH} -\sum_{i=1}^n (U_XrR_iU_X^*)\otimes T_i^*\right)^{-1}
(U_X\iota) \otimes I_\HH \\
&\qquad= (U_X\otimes I_{\DD_{T^*}}) K_{T,r}.
\end{split}
\end{equation*}
Thus
\begin{equation}\label{eq:assoc3}
\L_T\circ (\L_X\circ \L_\Rr)=
\begin{pmatrix}
 (U_X\otimes I_{\DD_{T^*}}) K_{T,r}&  (U_X\otimes I_{\DD_{T^*}}) \Theta_T(rR) (U_X^*\otimes I_{\DD_T})   \\
0 & 0
\end{pmatrix},
\end{equation}
and we have thus computed the right hand side of~\eqref{eq:assoc2}.

As for the left hand side, let us first remark that, computing $\L_T\circ \L_X$
according to~\eqref{eq:rprod}, we obtain as $(2,1)$ entry $I_\Fock\otimes(\alpha_{\L_{Y^*}}(T))^*$.
To avoid messy computations, we will use Lemma~\ref{le:structure} to obtain its other
entries.

Noting that in $\L_T$ and $\L_X$ all spaces have $\Fock$ as a tensor factor, 
and all operators have $I_\Fock$ as a factor, we shall write 
(a slight abuse of notation) $\L_T=I_\Fock\otimes \L'_T$, 
$\L_X=I_\Fock\otimes \L'_X$. Since both $\L'_T$ and $\L'_X$ are
unitary operators, the same is true of $\L'_T\circ \L'_X$.
Its $(2,1)$ entry is 
\[
c+dT^*(I-bT^*)^{-1}a=
\left(c^* +a^* T(I-b^* T)^{-1}d^*  \right)^*=
(\alpha_{\L_{Y^*}}(T))^*,
\]
while its $(1,1)$ entry is $D_{T^*}(I-bT^*)^{-1}a$. This last operator
has obviously dense range from $\HH$ to $\DD_{T^*}$ (remember that $a\not=0$),
and we may therefore apply Lemma~\ref{le:structure}. Consequently, the 
operators $Z_*:\DD_{\alpha_{\L_{Y^*}}(T)^*}\to \DD_{T^*}$ and
$Z:\DD_{\alpha_{\L_{Y^*}}(T)}\to \DD_T$, defined by
\[
\begin{split}
Z_*D_{\alpha_{\L_{Y^*}}(T)^*}&= D_{T^*}(I-bT^*)^{-1}a\\
Z D_{\alpha_{\L_{Y^*}}(T) }&= D_T(I-b^*T)d^*
\end{split}
\]
are unitary, and 
\[
\L'_T\circ \L'_X= 
\begin{pmatrix}
Z_*D_{\alpha_{\L_{Y^*}}(T)^*} & -Z_*\alpha_{\L_{Y^*}}(T) Z^*\\
\alpha_{\L_{Y^*}}(T)^* & D_{\alpha_{\L_{Y^*}}(T) }Z^*
\end{pmatrix}.
\]
Therefore
\[
\L_T\circ \L_X= 
\begin{pmatrix}
I_\Fock \otimes(Z_*D_{\alpha_{\L_{Y^*}}(T)^*} )& -I_\Fock \otimes (Z_*\alpha_{\L_{Y^*}}(T) Z^*)\\
I_\Fock \otimes \alpha_{\L_{Y^*}}(T)^* &I_\Fock \otimes (D_{\alpha_{\L_{Y^*}}(T) }Z^*)
\end{pmatrix}
=\L''\circ \L_{\alpha_{\L_{Y^*}}(T)},
\]
where $\L''=\begin{tpmatrix}I_\Fock\otimes Z_* & 0\\ 0 &I_\Fock\otimes Z^*\end{tpmatrix}$. Thus
\begin{align}\label{eq:assoc4}
(\L_T\circ \L_X)\circ \L_\Rr&=(\L''\circ \L_{\alpha_{\L_{Y^*}}(T)})\circ \L_\Rr
=\L''\circ (\L_{\alpha_{\L_{Y^*}}(T)}\circ \L_\Rr)\notag
\\
&=\L''\circ \begin{pmatrix} K_{\alpha_{\L_{Y^*}}(T),r} & \Theta_{\alpha_{\L_{Y^*}}(T)}(rR)\\
0 & 0
\end{pmatrix}\notag\\
&=
\begin{pmatrix}(I_\Fock\otimes Z_*) K_{\alpha_{\L_{Y^*}}(T),r} &(I_\Fock\otimes Z_*) 
 \Theta_{\alpha_{\L_{Y^*}}(T)}(rR) (I_\Fock\otimes Z^*)\\
0 & 0
\end{pmatrix}\\
&=
\begin{pmatrix}(I_\Fock\otimes Z_*) K_{\Phi_X^{-1}(T),r} &(I_\Fock\otimes Z_*) 
 \Theta_{\Phi_X^{-1}(T)}(rR) (I_\Fock\otimes Z^*)\\
0 & 0
\end{pmatrix}\notag
\end{align}
(we have used~\eqref{eq:inverse2} for the last equality.
We compare now~\eqref{eq:assoc3} with~\eqref{eq:assoc4}, and make $r\to1$. The resulting limits exist by~\eqref{eq:Kr} and~\eqref{eq:Thr}, and we obtain 
the assertion of the theorem.
\end{proof}

\section{Constrained row contractions}\label{se:constrained}

We introduce now some definitions from~\cite{GP3, GP4}, where the notion
of \emph{constrained} objects appears. Let $\Jfr$ be a WOT-closed two-sided
ideal in $\Lfr_n$, $\Jfr\not=\Lfr_n$.  We define two subspaces of $\Fock$
by 
\[
\MM_\Jfr= \overline{\Jfr\Fock}, \quad \NN_\Jfr= \Fock\ominus \MM_\Jfr.
\]
Then $\MM_\Jfr$ and $\Flip\MM_\Jfr$ are invariant to $L$ and to $R$, 
while $\NN_\Jfr$ and $\Flip\NN_\Jfr$
are invariant to $L^*$ and $R^*$.

The \emph{constrained left and right creation operators} belong to $\LL(\NN_\Jfr)$
and are given by 
\[
L_i^{\Jfr} = P_{\NN_\Jfr}L_i|\NN_\Jfr, \quad R_i^{\Jfr} = P_{\NN_\Jfr}R_i|\NN_\Jfr.
\]
An operator $M\in \LL(\NN_\Jfr\otimes \EE, \NN_\Jfr\otimes \EE_*)$ is called
\emph{multianalytic} if
\[
M(L_i^{\Jfr}\otimes I_\EE) = (L_i^{\Jfr}\otimes I_{\EE_*})M.
\]

We want to define  \emph{constrained row contractions} by using the
functional calculus with respect to elements of the ideal. A problem
appears, since for a general multicontraction the functional 
calculus is only defined for elements in $\lfr_n$; it can be
extended to $\Lfr$ only for completely noncoisometric contractions.
Thus, if $T$ is a general multicontraction, and $\mathfrak j\subset\lfr_n$
is a two-sided norm closed ideal, we say that $T$ is $\mathfrak j$-constrained
if $f(T)=0$ for all $f\in \mathfrak j$. If $T$ is c.n.c., and $\Jfr\subset\Lfr_n$,
we say that $T$ is $\Jfr$-constrained if $f(T)=0$ for all $f\in \Jfr$.
If $\Jfr$ is the wot-closure of $\mathfrak j$, and $T$ is c.n.c., then 
it is $\mathfrak j$-constrained iff it is $\Jfr$-constrained.

The next result connects the constraints with the automorphisms.

\begin{proposition}\label{pr:constr1}
If $T$ is $\mathfrak j$-constrained, then $\Phi_X^{-1}(T)$ is $\Phi_X(\mathfrak j)$-constrained (and similarly for $\Jfr$-constraints, in case $T$ is c.n.c.).
\end{proposition}

\begin{proof} 
If we denote $T'=\Phi_X^{-1}(T)$, and $\rho:\lfr\to \LL(\HH)$ is, as above,
the functional calculus for $T'$, then 
\[
T'_i=\rho((\Phi_{X}^{-1}(L))_i)=\rho(U_X^* L_i U_X)=\rho(\Phi_{X}^{-1}(L_i)).
\] 
Since
the functional calculus for $T'$ is the unique homomorphism algebra that
maps $L_i$ into $T'_i$, it must be $\rho\circ\Phi_{X}^{-1}$, and therefore
\[f(T')=\rho(\Phi_{X}^{-1}(f))=(\Phi_{X}^{-1}(f))(T) .\]
Thus $f(T')=0$ is equivalent to 
$(\Phi_{X}^{-1}(f))(T)=0$, whence the statement of the proposition follows.
\end{proof}


Now, if $T$ is a $\mathfrak j$-constrained contraction, and $\Jfr$ is the wot-closure of $\mathfrak j$, we define as in~\cite{GP3}:

\begin{enumerate}

\item[(a)] the \emph{constrained Poisson kernel}
$K_{\Jfr,T}:\HH\to\NN_\Jfr\otimes \DD_{T^*}$ by 
\[
K_{\Jfr,T}=P_{\NN_\Jfr\otimes\DD_{T^*}}K_T;
\]

\item[(b)]
the \emph{constrained characteristic function} 
$\Theta_{\Jfr,T}:\NN_\Jfr\otimes\DD_T \to \NN_\Jfr\otimes\DD_{T^*}$
by
\[
\Theta_{\Jfr,T}=P_{\NN_\Jfr\otimes\DD_{T^*}} \Theta_T |\NN_\Jfr\otimes\DD_T.
\]
\end{enumerate}

We can then give the following consequence of Theorem~\ref{th:theta} for
$\mathfrak j$-constrained multicontractions.

\begin{theorem}\label{th:theta2}
Suppose  $T$ is $\mathfrak j$-constrained, and denote $\mathfrak j'=\Phi_X(\mathfrak j)$, $\Jfr'=\Phi_X(\Jfr)$, $T'=\Phi_X^{-1}(T)$. Then
$\Theta_{\mathfrak J',T'}= (U_X\otimes \Om_*^*) \Theta_{\mathfrak J, T}(U_X^*\otimes \Om)$
and $K_{\Jfr',T'}=(U_X\otimes \Om_*^*) K_{\Jfr,T}$.
\end{theorem}

\begin{proof} 
By Theorem~\ref{th:theta}, we have
\[
\Theta_{T'}=( U_X \otimes \Om_*^* )\Theta_T ( U_X^*  \otimes \Om  ),
\]
where $\Om:\DD_{T'}\to \DD_T$ and 
$\Om_*:\DD_{T'{}^*}\to \DD_{T^*}$
are unitary operators.

We have 
\[
P_{\NN_{\mathfrak j'}}=U_X P_{\NN_{\mathfrak j}}U_X^*,\qquad
P_{\DD_{T'}}=\Om^* P_{\DD_T} \Om,\qquad
P_{\DD_{T'{}^*}} = \Om_*^* P_{\DD_{T^*}}\Om_*
\]
whence
\[
\begin{split}
P_{\NN_{\mathfrak j'}\otimes\DD_{T'} }&=(U_X\otimes\Om^* ) 
P_{\NN_{\mathfrak j}\otimes\DD_T}(U_X^*\lambda\otimes \Om),\\
P_{\NN_{\mathfrak j'}\otimes\DD_{T'{}^*} }&=(U_X\otimes \Om_*^*) 
P_{\NN_{\mathfrak j}\otimes\DD_{T^*}}(U_X^*\otimes \Om_*).
\end{split}
\]

Therefore
\[
\begin{split}
\Theta_{\mathfrak j',T'}
&= P_{\NN_{\mathfrak j'}\otimes\DD_{T'{}^*} }
\Theta_{T'} P_{\NN_{\mathfrak j'}\otimes\DD_{T'} }\\
&=(U_X\otimes \Om_*^*) 
P_{\NN_{\mathfrak j}\otimes\DD_{T^*}}(U_X^*\otimes \Om_*) 
\\&\qquad\qquad 
( U_X \otimes \Om_*^* )\Theta_T ( U_X^*  \otimes \Om  )
(U_X\otimes\Om^* ) 
P_{\NN_{\mathfrak j}\otimes\DD_T}(U_X^*\otimes \Om)\\
&=(U_X\otimes \Om_*^*) 
P_{\NN_{\mathfrak j}\otimes\DD_{T^*}}\Theta_T 
 P_{\NN_{\mathfrak j}\otimes\DD_T}(U_X^*\otimes \Om)\\
 &=(U_X\otimes \Om_*^*) \Theta_{\mathfrak j, T}(U_X^*\otimes \Om).
\end{split}
\]
The computations for the Poisson kernel are similar.
\end{proof}

If $\mathfrak j$ is the commutator ideal $\mathfrak c=[\lfr_n,\lfr_n]\subset \lfr_n$
(and correspondingly $\Cfr=[\Lfr_n,\Lfr_n]\subset\Lfr_n$, then the constraint
becomes just commutativity. Then $\Phi_X (\mathfrak c)=\mathfrak c$, $\Phi_X (\mathfrak C)=\mathfrak C$ for all~$X$, which translates in the fact that applying the automorphism $\Phi_X$ 
to a commuting multicontraction produces also a commuting multicontraction.
Theorem~\ref{th:theta2} yields then the transformation rule of the commutative characteristic
function with respect to automorphisms of the ball, as shown in~\cite{BT3}
(see Theorem 6.3 therein).

\section{Invariant ideals}

In connection to Theorem~\ref{th:theta2}, it is interesting to discuss
bilateral ideals $\Jfr$ of $\Lfr_n$ which are invariant with respect to 
all automorphisms $\Phi_X$. They have the property that, if $T$ is a $\Jfr$-constrained multicontraction,
$\alpha_X(T)$ is then also a $\Jfr$-constrained multicontraction.

We have already encountered the commutator ideal $\mathfrak C$. Other 
examples of invariant ideals are given by the
iterated
commutators $\mathfrak C^k$, defined by $\mathfrak C^{k+1}=
[\Lfr_n, \mathfrak C^k]$. These form
a decreasing sequence contained in $\mathfrak C$. We will prove
below that there are no invariant ideals
larger than $\mathfrak C$. But we need for this some more preparatory
results.

First, it is shown in~\cite{GP, GP1.5} that any $\Lfr_n$-invariant
subspace of $\Fock$ is of the form $\Theta (\Fock\otimes \EE)$, for 
$\EE$ a Hilbert space and $\Theta:\Fock\otimes \EE\to \Fock$ a multianalytic operator that
is also an isometry (such a $\Theta$ is called \emph{inner}). This multianalytic 
operator is essentially uniquely determined by the subspace: if 
$\Theta':\Fock\otimes \EE'\to \Fock$ satisfies 
$\Theta' (\Fock\otimes \EE')=\Theta (\Fock\otimes \EE)$, then there
exists a unitary $V:\EE'\to\EE$ such that $\Theta'=\Theta (I_\Fock\otimes V)$.
Based on these results, one proves
in~\cite{DP} that the map $\Jfr\mapsto\MM_\Jfr= \overline{\Jfr e_\emptyset} (=\overline{\Jfr\Fock})$
is a one to one map from the set of all bilateral ideals in
$\Lfr_n$ onto the set of subspaces in $\Fock$ invariant both to $\Lfr_n$
and to $\Rfr_n$. 

Finally, in~\cite{AP,DP2} one identifies the eigenvectors of $\Lfr_n^*$.
Namely, for any $\lambda\in\BBB^n$, one defines
\begin{equation}\label{eq:nu}
\nu_\lambda=(1-\|\lambda\|^2)^{1/2} (I-L[\bar\lambda])^{-1} e_\emptyset.
\end{equation}
Then $L_i^*\nu_\lambda=\bar\lambda_i \nu_\lambda$, whence
$\<L_w\nu_\lambda, \nu_\lambda\>=\lambda_w$ for any $w\in\FFF_n^+$.
Note that $\nu_\lambda$ are also eigenvectors of $\Rfr_n^*$ (corresponding to the same 
eigenvalues). The space spanned by all $\nu_\lambda$ ($\lambda\in\BBB^n$)
is $\MM_\Cfr^\perp$. This last space is the symmetric Fock space, which 
we will denote by $\Focks$, and 
the map $e_w\mapsto \lambda_w$ identifies it
with a space of functions on $\BBB^n$. If $A\in\Rfr_n$, then the projection 
of $A\1$ onto $\Focks$ is identified with the function $\<A\nu_\lambda, \nu_\lambda\>$.


\begin{theorem}\label{th:justc}
If $\mathfrak J\supset\Cfr$ is a bilateral ideal in $\Lfr_n$, and
$\Phi_X(\mathfrak J)=\mathfrak J$ for all $X$, then either 
$\mathfrak J=\mathfrak C$ or $\mathfrak J=\Lfr_n$.
\end{theorem}

\begin{proof}
Let $\MM=\overline{\Jfr \Fock}$ be the invariant subspace determined by $\Jfr$;
then $\Phi_X(\mathfrak J)=\mathfrak J$ implies $U_X\MM=\MM$. 
Suppose $\MM=\Theta (\Fock\otimes \GG)$;   define  
$\Gamma:\BBB^n\otimes\EE\to\CCC$ by the formula 
\begin{equation}\label{eq:Gamma}
\Gamma(\lambda,{h})=\<\Theta(\nu_\lambda\otimes {h}), \nu_\lambda\>.
\end{equation}
If $\Theta=\sum_{w} R_w\otimes m_w$, then 
\[
\begin{split}
\Gamma(\lambda,{h})&=\lim_{r\to1} \sum_w \<r^{|w|}(R_w\otimes m_w) (\nu_\lambda\otimes {h}), \nu_\lambda\> \\
&=\lim_{r\to1} \sum_w \<R_w\nu_\lambda, \nu_\lambda\>\<m_w {h},{e_\emptyset}\>=
\lim_{r\to1} \sum_w \lambda^w \<m_w {h},{e_\emptyset}\>.
\end{split}
\]
For $r<1$ the series on the right is uniformly convergent and thus defines
an analytic function  $\lambda\in\BBB^n$. It follows then that $\Gamma(\lambda, {h})$ is analytic
in $\lambda$. It is obviously linear in~${h}$; so we may consider $\lambda\mapsto \Gamma(\lambda,\cdot)$ as an analytic map $\tilde \Gamma$
from $\BBB^n$ into $\EE$ (actually, in the dual of $\EE$, which can be identified
with $\EE$).

On the other hand, since $U_X$ implements an automorphism of $\Rfr_n$,
one checks easily that $\Theta_X=U_X\Theta(U_X^*\otimes I_\EE)$ is also
an multianalytic inner operator. The invariance of $\MM$ with respect to 
$U_X$ implies that $\MM=\Theta_X(\Fock\otimes \EE)$. The essential uniqueness of this
representation implies then that for any $X\in U(1,n)$ there exists
$V_X\in\LL(\EE)$ such that
\[
U_X\Theta(U_X^*\otimes I_\EE)=\Theta(I_\Fock\otimes V_X).
\]

Let us take now $X$ such that $x>0$ (remember that the mappings $X\mapsto\Phi_X$ 
and $X\mapsto\phi_X$
have as kernel the constant unitaries). From~\eqref{eq:voic2} and~\eqref{eq:nu}
it follows then that 
\[
U_X^*e_\emptyset=U_X^*\nu_0=\nu_{\phi_X(0)}.
\]
Therefore
\[
\begin{split} 
\Gamma(\phi_X(0),{h})&=
\<\Theta (\nu_{\phi_X(0)}\otimes {h}), \nu_{\phi_X(0)}\>\\
&=\<U_X\Theta(U_X^*\otimes I_\EE)(e_\emptyset\otimes {h}), e_\emptyset\>\\
&=\<\Theta(I_\Fock\otimes V_X) (e_\emptyset\otimes {h}), e_\emptyset\>=
\Gamma(0,V_X({h})).
\end{split}
\]
This last relation can be rewritten as $\tilde \Gamma\circ \phi_X(0)=V_X^*\tilde \Gamma(0)$.
Since $V_X$ is unitary, we obtain that $\|\tilde \Gamma(\phi_X(0))\|=\|\tilde \Gamma(0)\|$.
The image of $\{X\in U(1,n):  x>0\}$ under the mapping $X\mapsto \phi_X(0)$ is the whole $\BBB^n$; therefore
$\tilde \Gamma$ is an analytic function on $\BBB^n$ with values in the Hilbert space $\EE$, of constant
norm, which must be actually constant.

For any ${h}\in\EE$ we can define an element $\Theta_{h}\in\Rfr$ by the formula
$\Theta_{h}\xi=\Theta(\xi\otimes h)$. We have then, by~\eqref{eq:Gamma} and the remarks before
the statement of the theorem, 
\[
\Gamma(\lambda, h)=\<\Theta_h\nu_\lambda, \nu_\lambda\>= (P_{\Focks}(\Theta_h{e_\emptyset}))(\lambda).
\]

Two cases present now. If $\tilde \Gamma$ is identically 0, then $(P_{\Focks}(\Theta_h{e_\emptyset}))=0$
for all $h\in\EE$, and thus the image of $\Theta$ is included in $\Cfr$. It follows that $\Jfr\subset\Cfr$;
and then the assumption implies $\Jfr=\Cfr$.

In the opposite case, take $h\in \EE$ such that $\Gamma(\lambda, h)$ is a nonnull constant.
Then $(P_{\Focks}(\Theta_h{e_\emptyset}))(\lambda)$ is a nonnull multiple of ${e_\emptyset}$. Thus
$\MM$ contains a vector of the form $a{e_\emptyset}+\xi_0$, with $a\not=0$ and $\xi_0\in\MM_\Cfr$.
But the assumption $\Cfr\subset\Jfr$ implies $\MM_\Cfr\subset \MM$; therefore ${e_\emptyset}\in\MM$.
Since $\MM$ is invariant, it follows that $\MM=\Fock$, whence $\Jfr=\Lfr$.
\end{proof}

\end{document}